\documentclass[11pt]{article}
\usepackage[utf8]{inputenc}
\usepackage{amsfonts}
\usepackage{amsmath} 
\usepackage{amssymb}
\usepackage{amsthm}
\usepackage{tikz}
\usepackage[margin=.85in]{geometry}
\usepackage{mathtools}
\usepackage{bbm}
\usepackage{xcolor}
\usetikzlibrary{matrix}
\usetikzlibrary{tikzmark,calc}
\newdimen\numht
\newdimen\numwd

\usepackage{biblatex}
\newtheorem{theorem}{Theorem}[section]
\newtheorem{lemma}[theorem]{Lemma}
\newtheorem{proposition}[theorem]{Proposition}
\newtheorem*{proposition*}{Proposition}
\newtheorem{corollary}[theorem]{Corollary}
\newtheorem*{theerem}{Theorem}
\newtheorem{definition}[theorem]{Definition}
\newtheorem{example}[theorem]{Example}
\newtheorem{remark}[theorem]{Remark}
\newtheorem{conj}[theorem]{Conjecture}

\newcommand{\bth}{\begin{theorem}}
\renewcommand{\eth}{\end{theorem}}
\newcommand{\bpr}{\begin{proposition}}
\newcommand{\epr}{\end{proposition}}
\newcommand{\bde}{\begin{definition}}
\newcommand{\ede}{\end{definition}}
\newcommand{\blem}{\begin{lemma}}
\newcommand{\elem}{\end{lemma}}
\newcommand{\bco}{\begin{corollary}}
\newcommand{\eco}{\end{corollary}}
\newcommand{\prove}{\begin{proof}}
\newcommand{\done}{\end{proof}}

\newcommand{\ite}{\begin{itemize}}
\newcommand{\mize}{\end{itemize}}

\newcommand{\ben}{\begin{enumerate}}
\newcommand{\een}{\end{enumerate}}

\renewcommand{\P}{\mathbb{P}}
\newcommand{\R}{\mathbb{R}}
\newcommand{\N}{\mathbb{N}}

\newcommand{\E}{\mathbb{E}}
\newcommand{\Var}{\mathrm{Var}}
\newcommand{\Cov}{\mathrm{Cov}}
\newcommand{\Tr}{\mathrm{Tr}}

\title{Fluctuations of Eigenvalues for Generalized Patterned Gaussian Random Matrices}
\author{ Frederick Rajasekaran\\\texttt{fredr@stanford.edu}\\Stanford University}
\date{\today\\[10pt]
	\begin{flushleft}
	                                       \\[5pt]
	\end{flushleft}}
\begin{document}
\maketitle

\abstract{In this work, we study a class of random matrices which interpolate between the Wigner matrix model and various types of patterned random matrices such as random Toeplitz, Hankel, and circulant matrices. The interpolation mechanism is through the correlations of the entries, and thus these interpolating models are highly inhomogeneous in their correlation structure. Historically, the study of random matrices has focused on homogeneous models, i.e., those with imposed structure (such as independence of the entries), as such restrictions significantly simplify computations related to the models. However, in this paper we demonstrate that for these interpolating inhomogenous models, the fluctuations of the linear eigenvalue statistics are approximately Gaussian. To handle the difficulties that come with inhomogeneity in the entries, we incorporate combinatorial arguments and recent tools from non-asymptotic random matrix theory.}

\section{Introduction}
In this paper, we study the fluctuations of linear statistics of eigenvalues for certain patterned random matrix models with general correlation structures among their entries. If $X_n$ is an $n \times n$ random matrix, its linear eigenvalue statistics are defined as
\begin{equation*}
    \sum_{i=1}^n f(\lambda_i)
\end{equation*}
where $\lambda_1, \dots, \lambda_n$ are the eigenvalues of $X_n$ and $f$ is a fixed test function. In particular, linear eigenvalue statistics are integrals of test functions against the empirical spectral distribution of the random matrix $X_n$, which is defined as $\mu_n(X_n) = \sum_{i=1}^{n} \delta_{\lambda_i}$.

The linear eigenvalue statistics of various random matrix models have been extensively studied in the literature. One is often interested in the fluctuations of linear statistics, and Gaussian central limit theorem results for various random matrix models have recently been proved. For fluctuations results for Wigner matrix models, see \cite{Chatterjee2009}, \cite{AndersonZeitouni2006}, \cite{Johansson1998}, \cite{Soshnikov2002}, \cite{SinaiSoshnikov1998}, \cite{Scherbina2011}, \cite{LytovaPastur2009}, and \cite{SosoeWong2013} and the references therein. For random Toeplitz matrix models, see \cite{Chatterjee2009} and \cite{LiuSunWang2012}, and for circulant matrices see \cite{AdhikariSaha2017} and \cite{AdhikariSaha2018}. For a non-Gaussian fluctuations result, see \cite{KumarMauryaSaha2022} for an example with odd monomial test functions and random Hankel matrices. The fluctuations of eigenvalues for symmetric circulant and reverse circulant matrices were studied in \cite{MauryaSaha2021} and \cite{AdhikariSaha2017}. 

To the best of our knowledge, the current linear statistics fluctuations results in the literature do not consider any models in which the matrix entries are allowed to have some general correlation structure. In this paper, we study the linear statistics of models that resemble the symmetric random Toeplitz, circulant, reverse circulant, symmetric circulant, and Hankel models, except we allow weaker correlations among certain entries. This work is a generalization of Section 4.3 of \cite{Chatterjee2009} and the related results in \cite{AdhikariSaha2017}.

The generalized random Toeplitz model is defined as follows, and is a more general version of the model studied in \cite{FriesenLowe2013}. The corresponding generalized circulant, reverse circulant, symmetric circulant, and Hankel models are defined in Section \ref{OtherModelsIntro}.

\bde
For $n \in \mathbb{N}$, a random matrix $X = X_n$ is called a generalized random Toeplitz matrix if it is symmetric and has independent diagonals up to symmetry. The entries of $X$ are assumed to satisfy $\E(X_{ij}) = 0$.
\ede

An generalized $n \times n$ Toeplitz matrix is determined by $n$ independent random vectors of dimensions varying from 1 through $n$, $\{\mathbf{X}_k\}_{k=0}^{n-1}$ with $\mathbf{X}_k \in \mathbb{R}^{n-k}$. Then the entries of the matrix are given by $X_{ij} = [\mathbf{X}_{|i-j|}]_{\min(i,j)}$. That is,
\begin{equation*}
    X = \begin{bmatrix}
        [\mathbf{X}_0]_1 & [\mathbf{X}_1]_1 & [\mathbf{X}_2]_1 & \cdots & [\mathbf{X}_{n-1}]_1 \\
        [\mathbf{X}_1]_1 & [\mathbf{X}_0]_2 & [\mathbf{X}_1]_2 & \cdots & [\mathbf{X}_{n-2}]_2 \\
        [\mathbf{X}_2]_1 & [\mathbf{X}_1]_2 & [\mathbf{X}_0]_3 & \cdots & [\mathbf{X}_{n-3}]_3 \\
        \vdots & \vdots & \vdots & \ddots & \vdots \\
        [\mathbf{X}_{n-1}]_1 & [\mathbf{X}_{n-2}]_2 & [\mathbf{X}_{n-3}]_3 & \cdots & [\mathbf{X}_0]_{n}
    \end{bmatrix}.
\end{equation*}
We refer to the vector $\mathbf{X}_k$ as the $k$th diagonal of the generalized Toeplitz matrix. In \cite{AdhikariBoseMauryaToeplitzCorrelated2023}, Adhikari, Bose, and Maurya study the limiting spectral distribution of random Toeplitz matrices with a correlation structure, but the model they consider in their paper is different from this one (see their Assumption 1). 

The most important example is the Gaussian case where the diagonal vectors are all jointly Gaussian vectors. In this case, the distribution of the entries of the matrix is completely determined by their covariances. For notation, we use
\begin{equation*}
c_k(i, j) = \Cov([\mathbf{X}_k]_i, [\mathbf{X}_k]_j) = \E([\mathbf{X}_k]_i[\mathbf{X}_k]_j).
\end{equation*}
Note that the values $c_k(i, j)$ also depend on $n$ and they cannot be chosen with complete freedom, since the $(n - k) \times (n - k)$ covariance matrix of each diagonal must be positive semidefinite. Also, due to the independent diagonals condition the covariance between any two entries on different diagonals is 0. The covariances between entries in the upper triangle of an Toeplitz matrix can thus be written as
\begin{equation*}
    \Cov(X_{ij}, X_{kl}) = \E([\mathbf{X}_{|i-j|}]_i [\mathbf{X}_{|k-l|}]_k) = \delta_{|i-j|=|k-l|}c_{|i-j|}(i, k), \qquad i \leq j,\: k \leq l.
\end{equation*}
In \cite{AjankiErdosKrugerGaussianCorrelated2016}, the authors consider another model for Gaussian matrices with a general correlation structure, though it is different from the one considered here.

As stated above, the covariances $c_k(i, j)$ completely determine the distribution of the entries in the Gaussian case. When $c_k(i, j) = \delta_{ij}$ for all $i, j, k$, we have a standard Gaussian Wigner matrix and when $c_k(i, j) = 1$ for all $i, j, k$, we have a standard Gaussian symmetric random Toeplitz matrix. Hence, one can think of generalized random Toeplitz matrices as a model that interpolates between the well-studied Wigner and random Toeplitz matrix models. 

In this paper, we prove Gaussian central limit theorems for the linear eigenvalue statistics with polynomial test functions of random Toeplitz matrices with generally correlated entries. We also prove a similar result for correlated versions of circulant, reverse circulant, symmetric circulant, and Hankel matrices (we give precise definitions of these models in Section \ref{OtherModelsIntro}). The convergence is with respect to the total variation distance between two random variables. Recall that for random variables $X$ and $Y$, 
\begin{equation*}
    d_{TV}(X, Y) = \sup_{A \in \mathcal{B}(\R)}|\P(X \in A) - \P(Y \in A)|,
\end{equation*}
where $\mathcal{B(\R)}$ is the set of Borel sets on $\R$. 

We then have the following theorem.

\bth
\label{MainTheorem}
Let $X_n$ be an $n \times n$ generalized random Toeplitz, circulant or symmetric circulant matrix with centered Gaussian entries with covariances labeled by $c_k(i, j)$. Assume that there exists constants $m$, $M \in (0, \infty)$ such that $m \leq \E(X_{ij}^2) \leq M$ for all $i, j,$ and $n$. Further, assume that there exists $\gamma > 0$ such that $c_k(i, j) \geq \gamma$ for all $i, j, k$, and $n$. Let $p$ be a positive integer and let $W_n = \mathrm{Tr}(X_n^p)$. Then as $ n \rightarrow \infty$,
\begin{equation*}
\frac{W_n - \E (W_n)}{\sqrt{\Var(W_n)}} \text{ converges in total variation to } N(0,1).
\end{equation*}

In the case when $X_n$ is a generalized random reverse circulant or Hankel matrix, the theorem still holds under the further assumption that $p$ is restricted to be an even positive integer.
\eth

\begin{remark}
     The condition on the parity of $p$ is necesary in the case of the Hankel matrices with generally correlated entries. In the standard Hankel matrix model, it was shown in \cite{KumarMauryaSaha2022} that the fluctuations of the linear eigenvalue statistics converge to a non-Gaussian limit for odd monomial test functions. However, it is not clear whether this condition is required in the reverse circulant case. 
\end{remark}

\begin{remark}
     Note that we can allow $p$ to grow with $n$ as $p = o(\log n/ \log \log n)$. Further, the above theorem holds for any fixed polynomial $f$ with non-negative coefficients. This simply follows from the proof of Theorem \ref{MainTheorem}, and the argument is spelled out in the proof of Theorem 4.5 of \cite{Chatterjee2009}.
\end{remark}

Furthermore, the generalized reverse circulant, circulant, and Hankel models (defined in Section \ref{OtherModelsIntro}) are not symmetric and hence do not necessarily have real eigenvalues. However, it still makes sense to talk about their eigenvalue statistics and the resulting fluctuations.

The above theorem allows for very general correlation structures among the diagonals. The covariances along each diagonal can fluctuate wildly with $n$, and as long as they stay bounded away from zero, the Gaussian central limit theorem still holds. The next theorem examines the fluctuations of the linear eigenvalue statistics in the regime where the covariances among the entries uniformly converge to zero.

\bth
\label{DecayTheorem}
Let $X_n$ be an $n \times n$ generalized random Toeplitz, circulant, reverse circulant, symmetric circulant, or Hankel matrix with centered Gaussian entries with covariances labeled by $c_k(i, j)$. Assume that there exists constants $m$, $M \in (0, \infty)$ such that $m \leq \E(X_{ij}^2) \leq M$ for all $i, j$ and for all $n$. Further, assume $c_k(i,j) = o(n^{-1/3})$ for all $i, j, k,$ and $n$ with $i \neq j$. Let $p$ be a positive integer and let $W_n = \mathrm{Tr}(X_n^p)$. Then as $ n \rightarrow \infty$,
\begin{equation*}
\frac{W_n - \E (W_n)}{\sqrt{\Var(W_n)}} \text{ converges in total variation to } N(0,1).
\end{equation*}
\eth

For results for related models with correlation decay, see \cite{AjankiErdosKrugerStabilityCorrelated2019} and \cite{ErdosKrugerSchroderCorrDecay2019}. However, in these papers they do not consider the fluctuations of linear eigenvalue statistics, and their correlation decay assumption is imposed through conditions on the multivariate cumulants of the entries.

In order to prove Theorem \ref{MainTheorem} and Theorem \ref{DecayTheorem}, we use the results of Chatterjee in \cite{Chatterjee2009} to bound the total variation distance between $W_n$ and a Gaussian with the same mean and variance. Doing so reduces the problem to computing lower bounds on $\Var(W_n)$ and upper bounds on $\E||X_n||$ and $||\Sigma||$, where $||\cdot||$ denotes the operator norm and $\Sigma$ is the $n^2 \times n^2$ covariance matrix of $X_n$. In order to lower bound $\Var(W_n)$, we expand the variance into a sum of covariances and use Wick's Theorem (stated in Section \ref{OutlineVarianceSection}) to show each term is positive. Then, a combinatorial argument shows that $\Var(W_n) \geq Kn$ for some constant $K$. Bounding $||\Sigma||$ is done easily with the Gershgorin Circle Theorem (stated in Section \ref{OutlineVarianceSection}). The challenging part of the proof is bounding $\E||X_n||$ due to the inhomogoneity of the entries of $X_n$. In order to prove the bound, we use recent work in non-asymptotic random matrix theory, and in particular we use matrix concentration inequalities. The idea is to construct $X_n$ as a series of deterministic matrices multiplied by i.i.d. standard Gaussian random variables, and then apply a matrix concentration inequality for matrices of this type.

The set up of the paper is as follows. We first prove, in detail, Theorem \ref{MainTheorem} in the case when $X_n$ is a random Toeplitz matrix with general correlation structure. Thus in Section \ref{OutlineVarianceSection} we give an outline of the generalized Toeplitz matrix proof and compute some necessary bounds on the variance. In Section \ref{IDIDBound} we compute a bound on the operator norm of an arbitrary generalized Toeplitz matrix and then use this bound to complete the proof of Theorem \ref{MainTheorem}. In Section \ref{OtherModelSection}, we comment on how the Toeplitz proof of Theorem \ref{MainTheorem} can be adapted for each of the other matrix models. In Section \ref{DecaySection}, we prove Theorem \ref{DecayTheorem}. Finally, in Section \ref{UniversalitySection} we conjecture about how Theorem \ref{MainTheorem} may be extended, both by removing the conditions on the covariances and by generalizing to sub-Gaussian entries. 

\subsection{Definitions of Other Generalized Patterned Matrix Models}
\label{OtherModelsIntro}

In this section we define other versions of patterned random matrix models that allow for general correlation structures among the entries. In \cite{AdhikariSaha2017}, Adhikari and Saha proved Gaussian fluctuations results for the standard versions of the following matrices and here we will prove similar results for their generalized versions. In all of the following models, for $k = 1, \dots, n$, let
\begin{equation*}
c_k(i, j) = \Cov([\mathbf{X}_k]_i, [\mathbf{X}_k]_j) = \E([\mathbf{X}_k]_i[\mathbf{X}_k]_j).
\end{equation*}

\bde
An $n \times n$ circulant matrix is defined as 
\begin{equation*}
    C_n = \begin{bmatrix}
        x_1 & x_2 & x_3 & \cdots & x_{n-1} & x_n \\
        x_{n} & x_1 & x_2 & \cdots & x_{n-2} & x_{n-1} \\
        x_{n-1} & x_n & x_1 &  \cdots & x_{n-3} & x_{n-2} \\
        \vdots & \vdots & \vdots & \ddots & \vdots & \vdots \\
        x_2 & x_3 & x_4 & \cdots & x_{n} & x_{1} \\
    \end{bmatrix}.
\end{equation*}
In other words, $(C_n)_{ij} = x_{(j-i + 1) \mod n}$. The $(j+1)$-th row is obtained by shifting the $j$-th row by one position to the right.
\ede

Note that in general, a circulant matrix is not symmetric and thus does not have real eigenvalues. However, the entries of the matrices will be restricted to be real, and we consider quantities of the form $\Tr(X_n^p)$ which will thus also be real. In order to introduce the correlations, we have $n$ independent centered random vectors $\mathbf{X}_k \in \R^n$ for $k = 1, \dots, n$. Then the generalized random circulant matrix can be written as
\begin{equation*}
    X_n = n^{-1/2} \begin{bmatrix}
        [\mathbf{X}_1]_1 & [\mathbf{X}_2]_1 & [\mathbf{X}_3]_1 & \cdots & [\mathbf{X}_{n-1}]_{1} & [\mathbf{X}_n]_1 \\
        [\mathbf{X}_{n}]_2 & [\mathbf{X}_1]_2 & [\mathbf{X}_2]_2 & \cdots & [\mathbf{X}_{n-2}]_{2} & [\mathbf{X}_{n-1}]_2 \\
        [\mathbf{X}_{n-1}]_3 & [\mathbf{X}_n]_3 & [\mathbf{X}_1]_3 &  \cdots & [\mathbf{X}_{n-3}]_{3} & [\mathbf{X}_{n-2}]_3 \\
        \vdots & \vdots & \vdots & \ddots & \vdots & \vdots \\
        [\mathbf{X}_2]_n & [\mathbf{X}_3]_n & [\mathbf{X}_4]_n & \cdots & [\mathbf{X}_{n-1}]_{n} & [\mathbf{X}_{1}]_n \\
    \end{bmatrix}.
\end{equation*}

Hence the generalized circulant matrix preserves the structure of the circulant matrix, but we allow for non-identical, correlated random variables along each diagonal. The next model is similarly defined.

\bde
An $n \times n$ reverse circulant matrix is defined as 
\begin{equation*}
    RC_n = \begin{bmatrix}
        x_1 & x_2 & x_3 & \cdots & x_{n-1} & x_n \\
        x_2 & x_3 & x_4 & \cdots & x_n & x_1 \\
        x_3 & x_4 & x_5 &  \cdots & x_1 & x_2 \\
        \vdots & \vdots & \vdots & \ddots & \vdots & \vdots \\
        x_n & x_1 & x_2 & \cdots & x_{n-2} & x_{n-1} \\
    \end{bmatrix}.
\end{equation*}
In other words, $(RC_n)_{ij} = x_{(i+j-1) \mod n}$. The $(j+1)$-th row is obtained by shifting the $j$-th row by one position to the left.
\ede

Note that while the reverse circulant matrix is similar the to the circulant matrix except rows are shifted to the left, the reverse circulant matrix is symmetric. However this property will break down in the generalized model. In order to introduce the correlations, we have $n$ independent centered random vectors $\mathbf{X}_k \in \R^n$ for $k = 1, \dots, n$. Then the generalized random reverse circulant matrix can be written as
\begin{equation*}
    X_n = n^{-1/2} \begin{bmatrix}
        [\mathbf{X}_1]_1 & [\mathbf{X}_2]_1 & [\mathbf{X}_3]_1 & \cdots & [\mathbf{X}_{n-1}]_{1} & [\mathbf{X}_n]_1 \\
        [\mathbf{X}_2]_2 & [\mathbf{X}_3]_2 & [\mathbf{X}_4]_2 & \cdots & [\mathbf{X}_n]_{2} & [\mathbf{X}_1]_2 \\
        [\mathbf{X}_3]_3 & [\mathbf{X}_4]_3 & [\mathbf{X}_5]_3 &  \cdots & [\mathbf{X}_1]_{3} & [\mathbf{X}_2]_3 \\
        \vdots & \vdots & \vdots & \ddots & \vdots & \vdots \\
        [\mathbf{X}_n]_n & [\mathbf{X}_1]_n & [\mathbf{X}_2]_n & \cdots & [\mathbf{X}_{n-2}]_{n} & [\mathbf{X}_{n-1}]_n \\
    \end{bmatrix}.
\end{equation*}

The generalized model is thus no longer symmetric unless it is maximally correlated. Next we have the symmetric circulant matrix in which the rows are shifted in the same manner as the circulant matrix, but we restrict the entries in each row to force the matrix to be symmetric. In this case, the generalized version will also stay symmetric.

\bde
An $n \times n$ symmetric circulant matrix is defined as 
\begin{equation*}
    SC_n = \begin{bmatrix}
        x_0 & x_1 & x_2 & \cdots & x_{2} & x_1 \\
        x_{1} & x_0 & x_1 & \cdots & x_{3} & x_{2} \\
        x_{2} & x_1 & x_0 &  \cdots & x_{4} & x_{3} \\
        \vdots & \vdots & \vdots & \ddots & \vdots & \vdots \\
        x_1 & x_2 & x_3 & \cdots & x_{1} & x_{0} \\
    \end{bmatrix}.
\end{equation*}
In other words, $(SC_n)_{ij} = x_{\frac{n}{2} - |\frac{n}{2}-|i-j||}$. The $(j+1)$-th row is obtained by shifting the $j$-th row by one position to the right.
\ede

In the generalized correlation model (similar to the standard random circulant model), the structure of the random matrix slightly depends on the parity of $n$. When $n$ is odd we have $(n+1)/2$ independent centered Gaussian random vectors $\mathbf{X}_k \in \R^n$ for $k = 0, \dots, (n-1)/2$. When $n$ is even, we have $n/2$ independent Gaussian random vectors $X_k \in \R^n$ for $k = 0, \dots \frac{n}{2} - 1$ and 1 independent Gaussian random vector $X_{n/2} \in R^{n/2}$. Then the generalized random symmetric circulant matrix can be written as
\begin{equation*}
    X_n = n^{-1/2} \begin{bmatrix}
        [\mathbf{X}_0]_1 & [\mathbf{X}_1]_1 & [\mathbf{X}_2]_1 & \cdots & [\mathbf{X}_{2}]_{n-1} & [\mathbf{X}_1]_n \\
        [\mathbf{X}_{1}]_1 & [\mathbf{X}_0]_2 & [\mathbf{X}_1]_2 & \cdots & [\mathbf{X}_{3}]_{n-1} & [\mathbf{X}_{2}]_n \\
        [\mathbf{X}_{2}]_1 & [\mathbf{X}_1]_2 & [\mathbf{X}_0]_3 &  \cdots & [\mathbf{X}_{4}]_{n-1} & [\mathbf{X}_{3}]_n \\
        \vdots & \vdots & \vdots & \ddots & \vdots & \vdots \\
        [\mathbf{X}_1]_n & [\mathbf{X}_2]_n & [\mathbf{X}_3]_n & \cdots & [\mathbf{X}_{1}]_{n-1} & [\mathbf{X}_{0}]_n \\
    \end{bmatrix}.
\end{equation*}

The following example exhibits the difference in structure due to the parity of the size of the matrix.

\begin{example}
We have 
\begin{equation*}
    X_4 = \frac{1}{2}\begin{bmatrix}
        [\mathbf{X}_0]_1 & [\mathbf{X}_1]_1 & [\mathbf{X}_2]_1 &  [\mathbf{X}_1]_4 \\
        [\mathbf{X}_{1}]_1 & [\mathbf{X}_0]_2 & [\mathbf{X}_1]_2 & [\mathbf{X}_{2}]_2 \\
        [\mathbf{X}_{2}]_1 & [\mathbf{X}_1]_2 & [\mathbf{X}_0]_3 &   [\mathbf{X}_{1}]_3 \\
        [\mathbf{X}_1]_4 & [\mathbf{X}_2]_2 & [\mathbf{X}_1]_3 & [\mathbf{X}_{0}]_4 \\
    \end{bmatrix} \quad \text{and} \quad
    X_5 = \frac{1}{\sqrt{5}}\begin{bmatrix}
        [\mathbf{X}_0]_1 & [\mathbf{X}_1]_1 & [\mathbf{X}_2]_1 &  [\mathbf{X}_{2}]_{4} & [\mathbf{X}_1]_5 \\
        [\mathbf{X}_{1}]_1 & [\mathbf{X}_0]_2 & [\mathbf{X}_1]_2 &  [\mathbf{X}_{2}]_{2} & [\mathbf{X}_{2}]_5 \\
        [\mathbf{X}_{2}]_1 & [\mathbf{X}_1]_2 & [\mathbf{X}_0]_3 &   [\mathbf{X}_{1}]_{3} & [\mathbf{X}_{2}]_3 \\
        [\mathbf{X}_{2}]_4 & [\mathbf{X}_2]_2 & [\mathbf{X}_1]_3 &   [\mathbf{X}_{0}]_{4} & [\mathbf{X}_{1}]_4 \\
        [\mathbf{X}_1]_5 & [\mathbf{X}_2]_5 & [\mathbf{X}_2]_3 &  [\mathbf{X}_{1}]_{4} & [\mathbf{X}_{0}]_5 \\
    \end{bmatrix}.
\end{equation*}
\end{example}

The last model we consider in this paper is the Hankel matrix. Similar to the reverse circulant model, even though the ungeneralized model is symmetric, this property of Hankel matrices is lost when we introduce the arbitrary correlations.

\bde
An $n \times n$ Hankel matrix is defined as 
\begin{equation*}
    H_n = \begin{bmatrix}
        x_1 & x_2 & x_3 & \cdots & x_{n-1} & x_n \\
        x_{2} & x_3 & x_4 & \cdots & x_{n} & x_{n+1} \\
        x_{3} & x_4 & x_5 &  \cdots & x_{n+1} & x_{n+2} \\
        \vdots & \vdots & \vdots & \ddots & \vdots & \vdots \\
        x_n & x_{n+1} & x_{n+2} & \cdots & x_{2n-2} & x_{2n-1} \\
    \end{bmatrix}.
\end{equation*}
In other words, $(H_n)_{ij} = x_{i+j-1}$. 
\ede

In order to introduce the correlations we have $n$ vectors $\mathbf{X}_k \in \R^k$ for $k = 1, \dots, n$ and $n-1$ vectors $\mathbf{X}_k \in \R^{2n-k}$ for $k = n+1, \dots 2n-1$. Then the generalized random Hankel matrix is

\begin{equation*}
 X_n = n^{-1/2} \begin{bmatrix}
        [\mathbf{X}_1]_1 & [\mathbf{X}_2]_1 & [\mathbf{X}_3]_1 & \cdots & [\mathbf{X}_{n-1}]_{1} & [\mathbf{X}_n]_1 \\
        [\mathbf{X}_{2}]_2 & [\mathbf{X}_3]_2 & [\mathbf{X}_4]_2 & \cdots & [\mathbf{X}_{n}]_{2} & [\mathbf{X}_{n+1}]_1 \\
        [\mathbf{X}_{3}]_3 & [\mathbf{X}_4]_3 & [\mathbf{X}_5]_3 &  \cdots & [\mathbf{X}_{n+1}]_{2} & [\mathbf{X}_{n+2}]_1 \\
        \vdots & \vdots & \vdots & \ddots & \vdots & \vdots \\
        [\mathbf{X}_n]_n & [\mathbf{X}_{n+1}]_{n-1} & [\mathbf{X}_{n+2}]_{n-2} & \cdots & [\mathbf{X}_{2n-2}]_{2} & [\mathbf{X}_{2n-1}]_1 \\
    \end{bmatrix}.
\end{equation*}

\subsection{Acknowledgements}
The author is grateful to Todd Kemp for suggesting this problem, supervising the project, and providing comments on the draft.

\section{Proof of Theorem 1.2 for Toeplitz Models Outline and First Steps}
\label{OutlineVarianceSection}

This section and the following section are devoted to proving Theorem \ref{MainTheorem} for generalized random Toeplitz matrices with correlated entries. We will prove the theorem under the assumption that all entries of $X_n$ are $N(0,1)$, since the proof is similar to the general variance case. We will comment where changes in the simplified proof can be made to accommodate general variances that are bounded below and above by $m$ and $M$, respectively. First, we include two important results from the literature that will be useful in our analysis.

\begin{theerem}[Gershgorin Circle Theorem]

Let $X$ be an $n \times n$ matrix with complex entries $x_{ij}$. For $i = 1, \dots, n$, let $R_i= \sum_{j \neq i}|x_{ij}|$. Then the spectrum of $X$ is contained in $\bigcup_{i=1}^{n} \overline{B_{R_i}(x_{ii})}$. In particular, for Hermitian matrices $||X||_{op} \leq \sup_{i} (|x_{ii}| + R_{i})$.
\end{theerem}

\begin{theerem}[Wick's Theorem]
    If $(X_1, \dots, X_n)$ is a mean-zero multivariate Gaussian random vector, then
    \begin{equation*}
        \E[X_1 X_2 \dots X_n] = \sum_{\pi \in P_2(n)}\prod_{\{i, j\} \in \pi}\E[X_i X_j],
    \end{equation*}
    where the sum is taken over all pair partitions of $\{1, \dots, n\}$.
\end{theerem}

We will prove the theorem assuming the entries of $X_n$ are standard Gaussian, but the generalization to general second moments bounded from above and below can easily be made, and in order to prove the theorem we use the technique of \cite{Chatterjee2009} in which Chatterjee gave a method for proving central limit theorems for linear statistics of eigenvalues of random matrices via second order Poincar\'e inequalities. The main result of \cite{Chatterjee2009} we will use is the following.

Let $(X_{ij})_{1 \leq i,j \leq n}$ be a collection of jointly Gaussian random variables with $n^2 \times n^2$ covariance matrix $\Sigma$. Let $X = n^{-1/2}(X_{ij})_{1\leq i, j \leq n}$. Then we have the following proposition.

\begin{proposition*}[Proposition 4.4. of \cite{Chatterjee2009}]
\label{ChatterjeeGaussian}
Take an entire function $f(z) = \sum_{m=0}^{\infty}b_m{z}^m$ and define $f_1$, $f_2$ as 
\begin{equation*}
f_1(z) = \sum_{m=1}^{\infty}m|b_m|z^{m-1} \qquad f_2(z) = \sum_{m=2}^{\infty}m(m-1)|b_m|z^{m-2}.
\end{equation*}
Let $\lambda$ denote the operator norm of $X$. Let $a = (\E f_1(\lambda)^4)^{1/4}$ and $b = (\E f_2(\lambda)^4)^{1/4}$. Suppose W = $\mathrm{Re Tr}f(X)$ has finite fourth moment and let $\sigma^2 = \Var(W)$. Let $Z$ be a normal random variable with the same mean and variance as $W$. Then 

\begin{equation}\label{ChatEqn}
d_{TV}(W,Z) \leq \frac{2\sqrt{5}||\Sigma||^{3/2}ab}{\sigma^2 n}.
\end{equation}
\end{proposition*}

We use this proposition to prove Theorem \ref{MainTheorem}. Here we consider monomials, so $f = x^p$, $f_1 = px^{p-1}$, and $f_2 = p(p-1)x^{p-2}$, and $X$ will be an $n \times n$ generalized random Toeplitz matrix (which we denote as $X_n$) with covariances in $[\gamma, 1]$. We need to bound three terms, $||\Sigma||$ (where $|| \cdot||$ is the operator norm), $ab$, and $\sigma^2$. The bound for $||\Sigma||$ is almost the same as in \cite{Chatterjee2009} for random Toeplitz matrices, and the argument for bounding $\sigma^2$ from below is also modeled off of the argument in \cite{Chatterjee2009} with some modifications to deal with the general covariances. The bound for $ab$ is new, and in order to bound the term we need to employ the matrix concentration inequalities in \cite{Tropp2018}. 

First, we give the proof of the bound on the operator norm of the covariance matrix $\Sigma$ via the Gershgorin circle theorem.

\blem \label{CovarGershgorinBound}
For any $n \times n$ generalized random Toeplitz matrix with $\Sigma$ as above, $||\Sigma|| \leq 2n$.
\elem

\prove
Let $\sigma_{ij,kl} = \Cov(X_{ij},X_{kl})$, so $\sigma_{ij,kl} = 0$ if $|i-j| \neq |k-l|$ and if $|i - j| = |k -l|$, $\sigma_{ij,kl} \leq 1$. By the Gershgorin circle theorem, if $\lambda_1 \dots \lambda_{n^2}$ are the eigenvalues of $\Sigma$, we have 
\begin{equation*}
||\Sigma|| = \max(|\lambda_m|) \leq \sup_{1\leq i, j \leq n} \sum_{k,l = 1}^n |\sigma_{ij,kl}|.
\end{equation*}
For any fixed $i,j$, there are at most $2n$ values of $k,l$ such that $|i - j| = |k - l|$. Hence each term in the maximum is bounded by $2n$ and thus so is $||\Sigma||$.
\done

\begin{remark}
    In the case where the entries of $X_n$ have general second moments, we have the bound $||\Sigma|| \leq 2Mn$.
\end{remark}

The proof for bounding the variance from below employs Wick's theorem to show that the covariances of products of the entries are positive, and then using the bound $\gamma$ to show the variance grows at least linearly.

\begin{remark}
This is the only part of the proof that uses the lower bound $\gamma$ of the covariances. This bound is sufficient to get linear growth of the variance, but may not be necessary. However, one cannot drop all conditions on the covariances and still maintain linear growth of the variance (see Corollary 1 or Theorem 2 of \cite{SinaiSoshnikov1998}, in which $\sigma^2$ converges to a finite limit for the Wigner matrix model).

%Does the fact that Sinai and Soshnikov have variance 1/4 matter at all? No, rescale by 2
\end{remark}
 
\blem
\label{VarBound}
With $W_n$ as in Theorem 1.2 for the Toeplitz model with general correlations, $\Var(W_n) \geq Kn$ for some constant $K$ that only depends on $p$ and $\gamma$.
\elem

\prove
It suffices to prove the lemma under the assumption that all of the entries are standard Gaussian random variables. In the case when they are not, normalize each entry, factor out the normalization coefficients, and apply the same argument (see the proof of Theorem 4.2 in \cite{Chatterjee2009}).

We first show that products of the entries all have nonnegative covariances. In this proof, $X_{ij} \coloneqq (X_n)_{ij}$ where $X_n$ is the generalized Toeplitz matrix. In the following computations, $n$ will be held fixed so we remove it from the notation. For any collections of non-negative integers $(\alpha_{ij})_{1 \leq i \leq j \leq n}$ and $(\beta_{ij})_{1\leq i \leq j \leq n}$ we have
\begin{equation}\label{pos_cov}
    \Cov \left(\prod X_{ij}^{\alpha_{ij}}, \prod X_{ij}^{\beta_{ij}} \right) = \E\left(\prod X_{ij}^{\alpha_{ij} + \beta_{ij}}\right) -  \E\left(\prod X_{ij}^{\alpha_{ij}}\right) \E\left(\prod X_{ij}^{\beta_{ij}}\right),
\end{equation}
and by the independent diagonals condition, these products factor as (with $ k = j - i$ and using the fact that $X_{ji} = X_{ij}$)
\begin{equation*}
\E\left(\prod X_{ij}^{\alpha_{ij} + \beta_{ij}}\right) -  \E\left(\prod X_{ij}^{\alpha_{ij}}\right) \E\left(\prod X_{ij}^{\beta_{ij}}\right) = \prod_{k = 0}^{n-1}\E\left(\prod X_{ij}^{\alpha_{ij} + \beta_{ij}}\right) - \prod_{k = 0}^{n-1}\E\left(\prod X_{ij}^{\alpha_{ij}}\right) \E\left(\prod X_{ij}^{\beta_{ij}}\right),
\end{equation*}
where the products on the right hand side inside the expectation are now taken over all $1 \leq i \leq j \leq n$ such that $j - i = k$. We show that for each $k$, the term in the product on the left is at least its corresponding term on the right. Fix $k$ (so all $i, j$ below are such that $j - i = k$) and consider
\begin{equation*}
\E\left(\prod X_{ij}^{\alpha_{ij} + \beta_{ij}}\right) - \E\left(\prod X_{ij}^{\alpha_{ij}}\right) \E\left(\prod X_{ij}^{\beta_{ij}}\right).
\end{equation*}
Then since the $X_{ij}$ are centered multivariate Gaussian random variables, by Wick's theorem if $\sum \alpha_{ij}$ or $\sum \beta_{ij}$ is odd, then the term on the right is 0. Thus it remains to consider the case when $\sum \alpha_{ij} = 2n$ and $\sum \beta_{ij} = 2m$ for some positive integers $m, n$. For notation purposes, enumerate the $2n+2m$ $X_{ij}$'s by $X_1, \dots X_{2n+2m}$ such that the first $X_1, \dots, X_{2n}$ correspond to an $X_{ij}$ which is raised to some $\alpha_{ij}$, and $X_{2n+1}, \dots X_{2n+2m}$ correspond to an $X_{ij}$ which is raised to some $\beta_{ij}$. Here, $P_2(k)$ denotes the set of pair partitions of $[k]$. In this case by Wick's Theorem
\begin{align*}
    \E\left(\prod X_{ij}^{\alpha_{ij} + \beta_{ij}}\right) - \E\left(\prod X_{ij}^{\alpha_{ij}}\right) \E\left(\prod X_{ij}^{\beta_{ij}}\right) 
    &= \sum_{\pi \in P_2(2(n+m))} \prod_{\{i, j\} \in \pi} \E(X_i X_j) \\
    & - \left(\sum_{\pi \in P_2(2n)} \prod_{\{i, j\} \in \pi} \E(X_i X_j) \right) \left(\sum_{\pi \in P_2(2m)} \prod_{\{i, j\} \in \pi} \E(X_{i+2n} X_{j+2n}) \right) .
\end{align*}
The map $\phi : P_2(2n) \times P_2(2m) \mapsto P_2(2(n+m))$ where $\phi(\pi, \sigma) = \pi \cup (2n +\sigma)$ (where addition is done element wise) is an injection. Thus every term in the double sum on the right has a corresponding term on the left. Since each of the $c_k(i, j) \geq 0$, $\E(X_i X_j) \geq 0$ for all $i, j$ so equation (\ref{pos_cov}) is nonnegative.

Now
\begin{equation*}
W_n = \mathrm{Tr}(X_n^p) = n^{-p/2} \sum_{1 \leq i_1, \dots, i_p \leq n} X_{i_1i_2}X_{i_2i_3}\dots X_{i_pi_1}.
\end{equation*}
Then since each of the terms above will have positive covariance, for any partition $\pi$ of any subset of $\{1, \dots, n\}^p$, 
\begin{align}
    \Var(W_n) &= n^{-p}\Var \left(\sum_{1 \leq i_1, \dots, i_p \leq n} X_{i_1i_2}X_{i_2i_3}\dots X_{i_pi_1}\right)  \nonumber \\
    &= n^{-p}\Var \left(\sum_{S \in \pi} \left(\sum_{(i_1, \dots, i_p) \in S} X_{i_1i_2}X_{i_2i_3}\dots X_{i_pi_1} \right) \right) \nonumber \\
    &\geq n^{-p} \sum_{S \in \pi} \Var \left(\sum_{(i_1, \dots, i_p) \in S} X_{i_1i_2}X_{i_2i_3}\dots X_{i_pi_1} \right). \label{var_final_ineq}
\end{align}

Now construct a partition by taking distinct positive integers $1 \leq a_1, a_2, \dots a_{p-1} \leq \lceil \frac{n}{3p} \rceil$ and let $D_{a_1, \dots, a_{p-1}}$ be the set of all $1 \leq i_1, \dots i_p \leq n$ such that $i_{k+1} - i_k = a_k$ for $k = 1, \dots, p-1$ and $1 \leq i_1 \leq \lceil \frac{n}{3} \rceil$. Then $|D_{a_1, \dots, a_{p-1}}| = \lceil \frac{n}{3} \rceil$ (since the choice of $i_1$ fixes all other $i_k$). Since the $a_i$'s are distinct, $X_{i_k i_{k+1}}$ is independent from $X_{i_j i_{j+1}}$ for all $j \neq k$. Thus
\begin{equation}\label{D_covar}
    \Var \left( \sum_{(i_1, \dots i_p) \in D} X_{i_1 i_2} \dots X_{i_p i_1}\right) = \sum_{(i_1, \dots i_p) \in D} \sum_{(i_1', \dots, i_p') \in D} \Cov(X_{i_1 i_2} \dots X_{i_p i_1}, X_{i_1' i_2'} \dots X_{i_p' i_1'}).
\end{equation}
Bounding the terms in the sum, using independence and the fact that the $X_{ij}$ are centered, 
\begin{align}
    \Cov(X_{i_1 i_2} \dots X_{i_p i_1}, X_{i_1' i_2'} \dots X_{i_p' i_1'}) &= \E(X_{i_1 i_2} \dots X_{i_p i_1} X_{i_1' i_2'} \dots X_{i_p' i_1'}) - \E(X_{i_1 i_2} \dots X_{i_p i_1}) \E(X_{i_1' i_2'} \dots X_{i_p' i_1'}) \nonumber \\
    &= \E(X_{i_1 i_2} X_{i_1' i_2'}) \dots \E(X_{i_p i_1} X_{i_p' i_1'}) \nonumber \\
    & \geq \gamma^p. \label{CovLowerBound}
\end{align}
Thus $(\ref{D_covar}) \geq \gamma^p |D|^2 \geq \frac{n^2 \gamma^p}{9}$.

The number of ways to choose $a_1, \dots, a_{p -1}$ satisfying the restrictions is 
\begin{equation*}
\lceil \frac{n}{3p} \rceil (\lceil \frac{n}{3p} \rceil - 1) \dots (\lceil \frac{n}{3p} \rceil - p + 2).
\end{equation*}
Since we can assume without loss of generality that $ n \geq 4p^2$, the above quantity can be lower bounded by $(n/12p)^{p-1}$ (since $\lceil \frac{n}{3p} \rceil - p + 2 \geq \frac{n}{3p} - p  \geq \frac{n}{3p} -\frac{n}{4p} = \frac{n}{12p}$). Finally note that if $(a_1, \dots, a_{p-1}) \neq (a_1', \dots a_{p-1}')$ then $D_{a_1, \dots, a_{p-1}}$ and $D_{a_1', \dots a_{p-1}'}$ are disjoint. Then applying (\ref{var_final_ineq}), 
\begin{equation*}
\Var(W_n) \geq n^{-p}\frac{n^{p-1}}{(12p)^{p-1}}\frac{n^2 \gamma^p}{9} = Kn.
\end{equation*}

\done

The last part of the proof is to get bounds on $ab$, which is done in the next section.

\section{Bounding the Spectral Norm of a Generalized Gaussian Toeplitz Matrix}\label{IDIDBound}

In order to bound the term $ab$ in (\ref{ChatEqn}), one must first bound the spectral norm of an arbitrary generalized random Toeplitz matrix. In order to bound the spectral norm, we write any generalized Gaussian random Toeplitz matrix as a sum of independent random matrices, and then apply a matrix concentration inequality to sums of this type. The following corollary of a matrix concentration inequality from \cite{Tropp2018} is used. If $H_1, \dots, H_k$ are fixed Hermitian matrices of common dimension $n$, and $\gamma_1, \dots, \gamma_k$ are standard normal random variables, then the random matrix
\begin{equation*}
X = \sum_{i = 1}^k \gamma_i H_i   
\end{equation*}
is called a \textit{Hermitian matrix Gaussian series}. The following result gives bounds on the operator norm of such matrices. 
\bth (Corollary 2.4 of \cite{Tropp2018}) \label{MatConcIneq}
Consider a Hermitian matrix Gaussian series $X = \sum_{i = 1}^k \gamma_i H_i$ with dimension $n$. Introduce the matrix standard deviation parameter
\begin{equation*}
\sigma (X) = ||\Var(X)||^{1/2} = \left \Vert\sum_{i = 1}^k H_i^2 \right \Vert^{1/2}.
\end{equation*}
Then
\begin{equation*}
    \frac{1}{\sqrt{2}} \cdot \sigma (X) \leq \E||X|| \leq \sqrt{e(1 + 2\log n)} \cdot \sigma (X),
\end{equation*}
where $||\cdot||$ denotes the spectral norm.
\eth

With this in mind, we can then get a bound on the spectral norm. 

\bth
\label{NormBoundIDID}
If $X_n$ is an $n \times n$ generalized random Toeplitz matrix and $\lambda_n$ is the spectral norm of $X_n$, then $\E(\lambda_n) \leq C\sqrt{\log n}$, where $C$ is a constant independent of $n$.
\eth

\begin{remark}
    Note that this theorem holds for any generalized Gaussian Toeplitz matrix, i.e. there are no restrictions on the covariances $c_k(i, j)$ as long as the corresponding covariance matrices are positive semidefinite.
\end{remark}

The first step in proving the theorem is to write any generalized Toeplitz matrix as a Hermitian matrix Gaussian series. Suppose $\mathbf{X}_k$ is a length $n-k$ vector of jointly Gaussian random variables with mean 0, variance 1, and covariance matrix $\Sigma$. If $A$ is a matrix such that $AA^T = \Sigma$ and $\mathbf{Z}_k$ is  length $n-k$ vector of i.i.d. standard normal random variables, then $A \mathbf{Z}_k$ has the same distribution as $\mathbf{X}_k$. Thus we can write
\begin{equation*}
    \mathbf{X}_k = \sum_{j = 1}^{n-k} [\mathbf{Z}_k]_j \mathbf{a}_{j}.
\end{equation*}
Here, $[\mathbf{Z}_k]_j$ are independent standard normal Gaussian random variables and $\mathbf{a}_{j}$ is the $j$th column of $A$. We will also denote $\mathbf{a}^{j}$ as the $j$th row of $A$.

\begin{remark}
\label{RemarkRowLength}
    Since each of the entries of the random matrix have variance at most $M$, this forces the Euclidean length of each row of $A$ to be at most $M$. This fact will be useful when bounding the operator norm of the matrix standard deviation parameter in the Gaussian series.
\end{remark}

\begin{remark}
    Such a decomposition of $\Sigma$ always exists since $\Sigma$ is positive semidefinite.
\end{remark}

Thus we can construct any diagonal of a generalized Toeplitz matrix in this manner, and constructing the full matrix is just a matter of piecing together the diagonals. For any diagonal vector $\mathbf{X}_k$ of a generalized Toeplitz matrix with Gaussian entries, let $\Sigma$ be its $(n-k) \times (n-k)$ covariance matrix and $A$ be such that $AA^T = \Sigma$. Then for $l = 1, \dots n-k$, let $B_{n,k,l}$ be an $n \times n$ matrix such that 
\begin{equation}
\label{BMatDef}
    (B_{n,k,l})_{ij} = \begin{cases}
        0 & |i - j| \neq k \\
        (A)_{min(i,j), l} & |i-j| = k.
    \end{cases}
\end{equation}
$B_{n, k, l}$ should be thought of as the $l$th column of $A$ pasted along the $k$th diagonal of an $n \times n$ matrix, with 0's elsewhere. The $A$ and $\Sigma$ matrices depend on $n$ and $k$ but for the sake of notation we omit the indices.
Then the random Toeplitz matrix with correlated entries can be written as
\begin{equation*}
    X_n = n^{-1/2}\sum_{k = 0}^{n-1}\sum_{l = 1}^{n-k} Z_{k,l}B_{n, k, l},
\end{equation*}
where each of the $Z_{k,l}$ are i.i.d. standard normal random variables.

\begin{example}
    For example, when $X_3$ is $3 \times 3$, the decomposition looks like the following. We remove the subscripts $n$ from the matrix entries for notational convenience. Let $A^{[k]}$ denote the square root matrix corresponding to the $k$th diagonal for $k = 0,1,2$.

\begin{align*}
    \sqrt{3} X_3 &= \sum_{l = 1}^3 Z_{0,l}\begin{bmatrix}
        (A^{[0]})_{1,l} & 0 & 0\\
         0 & (A^{[0]})_{2l} &  0\\
         0 & 0 & (A^{[0]})_{3l}\\
    \end{bmatrix}
    + \sum_{l=1}^{2}Z_{1,l}\begin{bmatrix}
         0 & (A^{[1]})_{1l} & 0\\
         (A^{[1]})_{1l} & 0 &  (A^{[1]})_{2l}\\
         0 & (A^{[1]})_{2l} & 0\\
    \end{bmatrix} \\
     &+ Z_{21}\begin{bmatrix}
        0 & 0 & (A^{[2]})_{11}\\
        0 & 0 & 0 \\
        (A^{[2]})_{11} & 0 & 0\\
    \end{bmatrix}
\end{align*}

Note that $(A^{[2]})_{11}$ is just the second moment of the corresponding entry in $\sqrt{3}X_3$.
\end{example}

We are now ready to prove Theorem 2.2.

\prove \textit{(Proof of Theorem \ref{NormBoundIDID})}

Again it suffices to prove the theorem under the assumption that the entries are standard Gaussians. The bounds for the general bounded variance case will just contain extra factors of $M$, which are independent of $n$.

In order to prove the theorem, we need to show that $||\sum_{k = 0}^{n-1}\sum_{l = 1}^{n-k}B_{n, k, l}^2||^{1/2} \lesssim n^{1/2}$. Since 
\begin{equation*}
    \left \Vert \sum_{k = 0}^{n-1}\sum_{l = 1}^{n-k}B_{n, k, l}^2 \right \Vert \leq \sum_{k = 0}^{n-1} \left \Vert \sum_{l = 1}^{n-k}B_{n, k, l}^2 \right \Vert,
\end{equation*}
it suffices to show that for any $k$, $||\sum_{l = 1}^{n-k}B_{n, k, l}^2|| \leq C$ for some constant C which is independent of $n$ and $k$ (in this case we will see that $C = 4$ works). The reason why such a bound holds is due to the fact stated in Remark \ref{RemarkRowLength}. If the $B_{n,k,l}$ were all diagonal matrices (as is the case when $k = 0$), then $\sum_{l = 1}^{k}B_{n, k, l}^2 = I_n$ which has operator norm 1. When the $B_{n,k,l}$ are not diagonal the entries ``mix" away from the main diagonal when the matrix is squared, and we need to bound the impact of this mixing on the operator norm of the matrix. 

To prove the above heuristics rigorously, we explicitly compute
\begin{equation}
\label{MatrixMult}
    (B_{n,k,l}^2)_{ij} = \sum_{p=1}^n (B_{n,k,l})_{ip}(B_{n,k,l})_{pj}.
\end{equation}
First consider the case when $0 < k \leq \frac{n}{2}$. Since $(B_{n,k,l})_{ij}$ is nonzero only when $|i-j| = k$, the term in the sum is nonzero when $|i-p| = |j-p| = k$. This condition can only hold for some $p$ if $|i -j| = 0$ or $|i-j|=2k$. Thus $(B_{n,k,l}^2)_{ij} = 0$ if $|i-j| \neq 0$ and $|i-j| \neq 2k$. Now consider the case when $i = j$. Again, here $A$ (which depends on $n$ and $k$) is the matrix associated to $B_{n, k, l}$, i.e. the $k$th diagonal of $B_{n, k, l}$ is the $l$th column of $A$. Then
\begin{align*}
    (B_{n,k,l}^2)_{ii} &= \sum_{p=1}^n (B_{n,k,l})_{ip}(B_{n,k,l})_{pi} \\
    &= (B_{n,k,l})_{i, i+k}(B_{n,k,l})_{i+k,i}\mathbbm{1}_{\{i \in [0, n-k]\}} + (B_{n,k,l})_{i, i-k}(B_{n,k,l})_{i-k, i}\mathbbm{1}_{\{i \in [k,n]\}} \\
    &= \begin{cases}
        ((A)_{il})^2 + ((A)_{i-k,l})^2 & i \in [k, n-k] \\
        ((A)_{il})^2 & i \in [0,k] \\
        ((A)_{i-k,l})^2 & i \in [n-k,n] .
    \end{cases}
\end{align*}
The last equality comes from equation (\ref{BMatDef}).

Now consider the case when $|i-j| = 2k$. Due to symmetry it suffices to consider the case in the upper triangle when $j = i+2k$. So $(B_{n,k,l}^2)_{i,i+2k} = \sum_{p=1}^n (B_{n,k,l})_{ip}(B_{n,k,l})_{p,i+2k}$, and for this to be nonzero requires $p = i + k$. Then we get 
\begin{align*}
    (B_{n,k,l}^2)_{i,i+2k} &= (B_{n,k,l})_{i, i+k}(B_{n,k,l})_{i+k, i+2k} \\
    &= (A)_{il}(A)_{i+k, l}.
\end{align*}
Putting this and the previous result together gives
\begin{equation*}
    (B_{n,k,l}^2)_{ij} = \begin{cases}
        ((A)_{il})^2 + ((A)_{i-k,l})^2 & i \in [k, n-k], i=j \\
        ((A)_{il})^2 & i \in [0,k], i = j \\
        ((A)_{i-k,l})^2 & i \in [n-k,n], i = j \\
        (A)_{\min(i,j),l}(A)_{\min(i,j)+k, l} & |i-j| = 2k \\
        0 & \text{otherwise}.
    \end{cases}
\end{equation*}
Then summing over $l$ and using Remark \ref{RemarkRowLength}, 
\begin{equation*}
    \left(\sum_{l=1}^{n-k} B_{n,k,l}^2 \right)_{ij} = \begin{cases}
        2 & i \in [k, n-k], i=j \\
        1 & i \in [0,k] \cup [n-k,n], i = j \\
        \mathbf{a}^{\min(i,j)} \cdot \mathbf{a}^{\min(i,j)+k} & |i-j| = 2k \\
        0 & \text{otherwise}.
    \end{cases}
\end{equation*}
By the Cauchy-Schwarz inequality and Remark \ref{RemarkRowLength}, $|\mathbf{a}^{\min(i,j)} \cdot \mathbf{a}^{\min(i,j)+k}| \leq 1$. Hence by the Gershgorin Circle Theorem, 
\begin{equation*}
\left \Vert \sum_{l = 1}^{n-k} B_{n,k,l}^2 \right \Vert \leq 4
\end{equation*}
since each Gershgorin disk is centered at 1 or 2 and has radius at most 2.

In the general case with non-unit variances, the above equation becomes
\begin{equation*}
    \left(\sum_{l=1}^{n-k} B_{n,k,l}^2 \right)_{ij} = \begin{cases}
        2M & i \in [k, n-k], i=j \\
        M & i \in [0,k] \cup [n-k,n], i = j \\
        \mathbf{a}^{\min(i,j)} \cdot \mathbf{a}^{\min(i,j)+k} & |i-j| = 2k \\
        0 & \text{otherwise}.
    \end{cases}
\end{equation*}
and $|\mathbf{a}^{\min(i,j)} \cdot \mathbf{a}^{\min(i,j)+k}| \leq M^2$ so 
\begin{equation*}
\left \Vert \sum_{l = 1}^{n-k} B_{n,k,l}^2 \right \Vert \leq 2M + 2M^2.
\end{equation*}

Back in the scenario with normalized entries, for the case when $k = 0$, each $B_{n, 0, l}$ is a diagonal matrix, so it follows from Remark 2.3 that $\sum_{l = 1}^n B_{n,0,l}^2 = I_n$ and thus $\left \Vert \sum_{l = 1}^n B_{n,0,l}^2 \right \Vert = 1$. When $k > \frac{n}{2}$, a similar computation to the case when $k \leq \frac{n}{2}$ can be done, and we see that there is no ``mixing" of matrix entries upon squaring. In this case, from equation (\ref{MatrixMult}), $(B_{n,k,l}^2)_{ij}$ is only nonzero when $|i-j| = 0$ or $|i-j| = 2k$. Since $k > \frac{n}{2}$, $|i-j| \neq 2k$ for any $i, j$, so $B_{n,k,l}^2$ is diagonal. When $i = j$ we have $(B_{n,k,l}^2)_{ii} = \sum_{p=1}^n (B_{n,k,l})_{ip}(B_{n,k,l})_{pi}$, and for the terms to be nonzero we need $p = i+k$ or $p = i-k$. $p = i+k$ can only happen when $i \in [0, n-k]$ and $p = i-k$ can only happen if $i \in [k,n]$. Thus, from equation (\ref{BMatDef}), 
\begin{equation*}
    (B_{n,k,l}^2)_{ij} = \begin{cases}
        ((A)_{il})^2 & i \in [0, n-k], i=j \\
        ((A)_{i-k, l})2 & \in i \in [k, n], i=j \\
        0 & \text{otherwise}.
    \end{cases}
\end{equation*}
Summing over $l$ and applying Remark \ref{RemarkRowLength}, we see that $\sum_{l = 1}^{n-k} B_{n,k,l}^2$ is the identity matrix with some diagonal elements replaced by 0's, and thus $\left \Vert \sum_{l = 1}^{n-k} B_{n,k,l}^2 \right \Vert \leq 1$.

Hence we have shown for any $k$, $||\sum_{l = 1}^{k}B_{n, k, l}^2|| \leq 4$, and thus $||\sum_{k = 0}^{n-1}\sum_{l = 1}^{n-k}B_{n, k, l}^2||^{1/2} \leq 2\sqrt{n}$ as desired.

Then, applying the Matrix Khintchine inequality, 
$$
    \E||X_n|| \leq n^{-1/2} \sqrt{e(1+2 \log n)}(2\sqrt{n}) \leq C\sqrt{\log n}
$$
for $n$ large enough and some constant $C$ which does not depend on $n$.
\done

\blem
\label{OperatorNormMomentBound}
$\E(\lambda_n^k) \leq (Ck\log n)^{k/2}$ for any $n$ and $k$, where $C$ is universal constant. 
\elem
\prove
This is proved via concentration of measure techniques (see \cite{Ledoux2001}). From equation (11.2) in \cite{Kemp2013}, 
\begin{equation*}
\mu(\{|\lambda_n - \E(\lambda_n)| \geq t\}) \leq 2e^{-t^2/2||F||_{Lip}^2},
\end{equation*}
where $||F||_{Lip}$ is the value of the Lipschitz function that gives the spectral norm of $A_n$ (which is a finite constant). Now we use the layer cake representation. Let $\kappa_k(dt) = kt^{p-1}dt$ on $[0, \infty)$. The corresponding cumulative function is $\phi_k(x) = \int_0^x d\kappa_k = x^k$. Then by Proposition 12.5 of \cite{Kemp2013} applied to the random variable $X_n = |\lambda_n - \E(\lambda_n)|$,
\begin{equation*}
\E(X_n^k) = \int_0^{\infty} \P(X_n \geq t)kt^{k-1}dt \leq 2k \int_0^{\infty}t^{k-1}e^{-t^2/||F||_{Lip}^2},
\end{equation*}
and then substituting $s = \frac{t}{\sqrt{2}||F||_{Lip}}$ the above becomes (with $C$ possibly changing between lines)
\begin{align*}
    2k\left(\sqrt{2}||F||_{Lip}^{k-1}\int_0^{\infty}s^{k-1}e^{-{s^2}}ds (\sqrt{2}||F||_{Lip})\right) &= 2k(2||F||_{Lip}^2)^{k/2} \int_0^{\infty}s^{k-1}e^{-s^2}ds\\
    &= C^{k/2}k\int_0^{\infty}s^{k-1}e^{-s^2}ds \\
    &= C^{k/2} k \Gamma\left(\frac{k}{2}\right) \\
    &\leq C^{k/2}k\left(\frac{k}{2}\right)^{k/2-1} \\
    &= C^{k/2}k^{k/2}.
\end{align*}

Then (again with $C$ possibly changing between inequalities)
\begin{align*}
    ||\lambda_n||_{k} &\leq ||\lambda_n - \E(\lambda_n)||_{k} + ||\E(\lambda_n)||_{k} \\ 
     &\leq C(k^{1/2} + \sqrt{\log n}) \\
     &\leq C(k^{1/2}\sqrt{\log n}),    
\end{align*}
where $||\cdot||_{k}$ denotes the $L^k$-norm, and the last inequality is for $n$ larger than $e^2$ and $k \geq 2$ (the case when $k = 1$ is Theorem \ref{NormBoundIDID}). 
\done

Now we can finally prove Theorem \ref{MainTheorem}.
\prove \textit{(Proof of Theorem \ref{MainTheorem})}

Following Proposition \ref{ChatterjeeGaussian}, $f_1(x) = px^{p-1}$ and $f_2(x) = p(p-1)x^{p-2}$, and applying H\"older's inequality gives $ab \leq p^3(\E(\lambda_n^{4p}))^{1/2}$.
Let $Z_n$ be a Gaussian random variable with the same mean and variance as $W_n$. Then by Proposition \ref{ChatterjeeGaussian} and Lemma \ref{CovarGershgorinBound},
\begin{equation*}
    d_{TV}(W_n, Z_n) \leq \frac{Cp^3 (\E(\lambda_n^{4p}))^{1/2}\sqrt{n}}{\Var(W_n)}
\end{equation*}
where $C$ is a universal constant. From Lemma \ref{OperatorNormMomentBound}, the $p^3(\E(\lambda_n^{4p}))^{1/2}$ term is bounded by $p^3(Cp \log n )^p$. Incorporating Lemma \ref{VarBound}, 
\begin{equation*}
d_{TV}(W_n, Z_n) \leq \frac{C^p p^{p+3}(\log n)^p}{\sqrt{n}}
\end{equation*}
and this goes to 0 as $n \rightarrow \infty$. Note that in this final step, the constant $C$ now depends on $p$ (it contains a factor of $\gamma^{-p}$). Centering and normalizing proves the theorem.
\done

\section{Fluctuations for Other Matrix Models with Correlated Entries}
\label{OtherModelSection}
The arguments and methods used above to prove the central limit theorem for generalized Toeplitz matrices extend to other matrix models. In \cite{AdhikariSaha2017}, Adhikari and Saha used Chatterjee's total variation bound to prove fluctuations results for circulant, symmetric circulant, reverse circulant, and Hankel matrices. In this section, we extend these results to the four corresponding matrix models that allow for general correlation structures among the entries. The proofs for each of these follow the same structure by computing upper bounds on $||\Sigma||$ and $\E||X_n||$ and lower bounds on $\Var(W_n)$. The bound for $\Sigma$ is essentially the same for all of the models. Furthermore, in this section we will assume that all entries of the random matrices have unit variance. The proofs can be adapted to the general second moments case in the same way as for the Toeplitz models.

\blem
\label{CovarianceMatrixBoundOtherModels}
Let $X_n$ be an $n \times n$ circulant, reverse circulant, symmetric circulant, or Hankel matrix with standard Gaussian correlated entries. Let $\Sigma$ denote its covariance matrix. Then $||\Sigma|| \leq 2n$.
\elem

\prove
For any $i, j$, $(X_n)_{ij}$ is correlated with at most $2n-1$ other entries of $X_n$. Then by the Gershgorin circle theorem, 
\begin{equation*}
    ||\Sigma|| \leq 1 + (2n-1) \sup_{k, i \neq j} c_k(i, j) \leq 2n.
\end{equation*}
\done

\subsection{Reverse Circulant Matrices with Correlated Entries}

Since the reverse circulant case is the longest proof of the four models, we prove it first. As stated above, we use Proposition \ref{ChatterjeeGaussian} and bound $\Var(X_n)$ from below and $\E[||X_n||]$ from above. 

\blem
With $W_n$ as in Theorem \ref{MainTheorem} and $X_n$ a reverse circulant matrix with correlated entries, $\Var(W_n) \geq Kn$ for some constant $K$ that only depends on $p$ and $\gamma$.
\elem
The following proof is a version of the proof of Lemma 10 in \cite{AdhikariSaha2017}. We can use the same general partition argument, but incorporate the bound $\gamma$ on the covariances. 

\prove
Again letting $X_{ij} \coloneqq (X_n)_{ij}$, we have 
\begin{equation*}
W_n = \Tr(X_n^p) = n^{-p/2} \sum_{1 \leq i_1, \dots, i_p \leq n} X_{i_1i_2}X_{i_2i_3}\dots X_{i_pi_1}.
\end{equation*}

Showing that all of the terms in the above sum are positively correlated is done in a similar way to the proof of Lemma \ref{VarBound}, except the products on the right hand side of equation (\ref{pos_cov}) factor with the condition of $(j+i-1) \mod n$ instead of $j-i$. Then the same argument using Wick's theorem shows that the terms in the above sum are positively correlated.

Then for any partition $\pi$ of any subset of $\{1, \dots, \frac{n}{3}\}^p$,
\begin{equation}   
\Var(W_n) \geq n^{-p} \sum_{S \in \pi} \Var \left(\sum_{(i_1, \dots i_p) \in S} X_{i_1i_2}X_{i_2i_3}\dots X_{i_pi_1} \right). \label{PartitionRC}
\end{equation}

Now we construct an appropriate partition following the proof of Lemma 10 in \cite{AdhikariSaha2017}. Let 
\begin{equation*}
\mathcal{A} = \left\{(a_1, \dots, a_p) \in \N^p : \frac{kn}{3p} + 1 \leq a_k \leq \frac{(k+1)n}{3p}, k = 1, 2, \dots, p-1 \right\}
\end{equation*}
and define
\begin{equation*}
D_{a_1, \dots, a_p} = \{(i_1, \dots, i_p) : 1 \leq i_1 \leq \frac{n}{3p}, i_k +i_{k+1} - 1 = a_k, k = 1, \dots, p\},
\end{equation*}
where $i_{p+1} = i_1$ and $(a_1, \dots, a_p) \in \mathcal{A}$. Note that we did not impose a condition on $a_p$ in the definition of $\mathcal{A}$. However, in the case when $p$ is even, we will show that given $(a_1, \dots, a_{p-1})$, the value of $a_p$ is fixed and uniquely determined by $(a_1, \dots a_{p-1})$. Also note that the $a_i$'s are distinct and because of this $X_{i_k i_{k+1}}$ is independent from $X_{i_j i_{j+1}}$ for $j \neq k$. Thus given some $D_{a_1, \dots, a_p}$ (which we denote as $D$ in the following computation),
\begin{equation*}
    \Var \left( \sum_{(i_1, \dots i_p) \in D} X_{i_1 i_2} \dots X_{i_p i_1}\right) = \sum_{(i_1, \dots i_p) \in D} \sum_{(i_1', \dots i_p') \in D} \Cov(X_{i_1 i_2} \dots X_{i_p i_1}, X_{i_1' i_2'} \dots X_{i_p' i_1'}).
\end{equation*}
Then following equation (\ref{CovLowerBound}), 
\begin{equation*}
    \Cov(X_{i_1 i_2} \dots X_{i_p i_1}, X_{i_1' i_2'} \dots X_{i_p' i_1'}) \geq \gamma^p,
\end{equation*}
so
\begin{equation*}
\Var \left( \sum_{(i_1, \dots i_p) \in D} X_{i_1 i_2} \dots X_{i_p i_1}\right) \geq \gamma^p |D|^2.
\end{equation*}

Note that if $(a_1, \dots, a_{p}) \neq (a_1', \dots a_{p}')$ then $D_{a_1, \dots, a_{p-1}}$ and $D_{a_1', \dots a_{p-1}'}$ are disjoint. Then from equation (\ref{PartitionRC}),

\begin{align*}
    \Var(W_n) &\geq n^{-p}\sum_{(a_1, \dots, a_p) \in \mathcal{A}} \Var\left(\sum_{(i_1, \dots, i_p) \in D} X_{i_1 i_2} \dots X_{i_p i_1}\right) \\
    &\geq n^{p}|\mathcal{A}||D|^2 \gamma^p. \\
\end{align*}

 It then remains to compute the sizes of the sets $\mathcal{A}$ and $D_{a_1, \dots, a_p}$. For a given set of values $(a_1, \dots, a_p)$, recursively applying the relation $i_k + i_{k+1} - 1 = a_k$ for $k = 1, \dots, p-1$ we get the following equation for $i_p$:
\begin{equation*}
i_p = \begin{cases}
    a_{p-1} - a_{p-2} + \dots - a_2 + a_1 - i_1 + 1 & \text{if $p$ is even,} \\
    a_{p-1} - a_{p-2} + \dots + a_2 -a_1 + i_1 & \text{if $p$ is odd.}
\end{cases}
\end{equation*}
Then applying the relation for $k = p$ gives $a_p = i_p + i_1 - 1 = a_{p-1} - a_{p-2} + \dots - a_2 + a_1$, when $p$ is even. Thus if $p$ is even then $a_p$ is determined by $a_1, \dots a_{p-1}$ and it does not depend on $i_1$, but if $p$ is odd then $a_p$ depends on $a_1, \dots, a_2, \dots a_{p-1}$ and $i_1$. Hence in the case when $p$ is even one may think of $\mathcal{A}$ as a collection of $(p-1)$-tuples $(a_1, \dots, a_{p-1})$ subject to the constraints in the definition of $\mathcal{A}$ above. We can then easily calculate that $|\mathcal{A}| = (\frac{n}{3p})^{p-1}$ since the values of $a_1, \dots, a_{p-1}$ are chosen from intervals of length $\frac{n}{3p}$.

Now, when $p$ is even, for a fixed choice of $a_1, \dots a_{p-1}$ the values of $i_2, \dots, i_p$ are uniquely determined by the choice of $i_1$, so the number of elements in $D_{a_1, \dots a_p}$ is the same as the number of ways of choosing $i_1$, so 
\begin{equation*}
|D_{a_1, \dots a_p}| = \frac{n}{3p}.
\end{equation*}
Thus, when $p$ is even
\begin{equation*}
\Var(W_n) \geq n^{-p} \left(\frac{n}{3p}\right)^{p-1} \left(\frac{n}{3p}\right)^2 = Kn
\end{equation*}
where $K$ only depends on $\gamma$ and $p$.
\done

The last part of the proof of Theorem \ref{MainTheorem} for reverse circulant matrices is a bound on the operator norm.

\blem
\label{NormBoundRC}
If $X_n$ is an $n \times n$ reverse circulant matrix with correlated entries and $\lambda_n$ is the spectral norm of $X_n$, then $\E(\lambda_n) \leq C\sqrt{\log n}$, where $C$ is a constant independent of $n$. 
\elem

In order to prove the above Lemma, we need to be slightly more careful than in the proof of Theorem \ref{NormBoundIDID}. Note that Theorem \ref{MatConcIneq} requires the matrices in the Gaussian series to be self-adjoint. However, generalized reverse circulant matrices with correlated entries are not symmetric, and thus neither will be the coefficient matrices from the construction in the beginning of Section \ref{IDIDBound}. 

We can get around this via the following dilation argument. If $A$ is an $n \times n$ matrix (not necessarily self-adjoint), then its \textit{dilation} is
\begin{equation*}
    \Tilde{A} \coloneqq \begin{bmatrix}
        0 & A^* \\
        A & 0 \\
    \end{bmatrix}.
\end{equation*}

Here, we are only considering real matrices, so $A^* = A^T$. Then since $\Tilde{A}$ is a $2n \times 2n$ self-adjoint matrix we have 
\begin{equation*}
    ||\Tilde{A}||^2 = ||\Tilde{A}^2|| = \max\{||A^T A||, ||AA^T||\} = ||A||^2
\end{equation*}
where $||\cdot||$ denotes the operator norm. Hence we will use the dilations of the coefficient matrices in the matrix series to compute the bound on the operator norm via matrix concentration inequalities. 

\prove
We construct the matrices $B_{n, k, l}$ in a similar manner as done in section 3. Recall that $n$ is the size of the matrix, $k$ labels the Gaussian vector $\mathbf{X}_k$, and $l$ corresponds to the $l$th column of $A$, which is such that $AA^T = \Sigma$ where $\Sigma$ is the covariance matrix of $\mathbf{X}_k$. Again, $A$ depends on $k$ and $n$, but we forgo the labeling for notation purposes. We then have 
\begin{equation}\label{RCBMatrix}
    (B_{n, k, l})_{ij} = \begin{cases}
        0 & i + j - 1 \not \equiv k \mod n \\
        A_{i l} & i + j -1 \equiv k \mod n \\
    \end{cases}.
\end{equation}
Then we can write
\begin{equation*}
    X_n = n^{-1/2}\sum_{k=1}^n \sum_{l=1}^n Z_{k,l}B_{n,k,l},
\end{equation*}
where the $Z_{k, l}$ are standard normal random variables.  Then applying the dilations, 
\begin{equation*}
    \Tilde{X_n} = n^{-1/2}\sum_{k=1}^n \sum_{l=1}^n Z_{k,l}\widetilde{B_{n,k,l}}.
\end{equation*}
Hence we must now show $\left \Vert \sum_{k=1}^n \sum_{l=1}^n\widetilde{B_{n,k,l}}^2 \right \Vert^{1/2} \lesssim n^{1/2}$. Now
\begin{equation*}
    \widetilde{B_{n, k, l}^2} = \begin{bmatrix}
        B_{n, k, l}B_{n, k, l}^T & 0 \\
        0 & B_{n, k, l}^TB_{n, k, l}
        \end{bmatrix}
\end{equation*}
so we compute
\begin{equation*}
    (B_{n, k, l} B_{n, k, l}^T)_{ij} = \sum_{i = 1}^n (B_{n, k, l})_{ip} (B_{n, k, l})_{jp}
\end{equation*}
and for this to be nonzero, from equation (\ref{RCBMatrix}) we need $i + p -1 \equiv k \mod n$ and $j + p - 1 \equiv k \mod n$ which forces $i = j$, so $B_{n, k, l} B_{n, k, l}^T$ and $B_{n, k, l}^T B_{n, k, l}$ are diagonal. Then from equation (\ref{RCBMatrix}),
\begin{equation*}
    (B_{n, k, l} B_{n, k, l}^T)_{ii} = A_{il}^2
\end{equation*}
and similarly for $B_{n, k, l}^T B_{n, k, l}$. Then by  the fact that the rows of $A$ have length 1, $\sum_{l=1}^n\widetilde{B_{n,k,l}}^2 = I_{2n \times 2n}$, and it follows that $\left \Vert \sum_{k=1}^n \sum_{l=1}^n\widetilde{B_{n,k,l}}^2 \right \Vert^{1/2} \leq n^{1/2}$ (in fact we have equality here).
\done

\subsection{Circulant Matrices with Correlated Entries}

The proof of the theorem again involves two parts along with Lemma \ref{CovarianceMatrixBoundOtherModels}.

\blem
With $W_n$ as in Theorem \ref{MainTheorem} and 
$X_n$ a circulant matrix with correlated entries, $\Var(W_n) \geq Kn$ for some constant $K$ that only depends on $p$ and $\gamma$.
\elem

\prove

This can be proved in similar manner to Lemma \ref{VarBound}. Notice that the upper triangle of a circulant matrix with correlated entries is identical to that of a Toeplitz matrix with correlated entries. The argument used in Lemma \ref{VarBound} only involves a partition involving the upper triangle (since $a_1, \dots, a_{p-1}$ are all positive and $i_{k+1} - i_{k} = a_k$), and the only point where the lower triangle is involved is the term $X_{i_{p} i_{1}}$. However, the bound in equation (\ref{CovLowerBound}) still holds (up to a constant factor) even if $X_{i_{p} i_{1}}$ is correlated with one other $X_{i_{k} i_{k+1}}$. To see this, suppose $X_{i_{p} i_{1}}$ is not independent of $X_{i_{k} i_{k+1}}$ for some $k$. Then
\begin{align*}
    \Cov(X_{i_1 i_2} \dots X_{i_p i_1}, X_{i_1' i_2'} \dots X_{i_p' i_1'}) &= \E(X_{i_1 i_2} \dots X_{i_p i_1} X_{i_1' i_2'} \dots X_{i_p' i_1'}) - \E(X_{i_1 i_2} \dots X_{i_p i_1}) \E(X_{i_1' i_2'} \dots X_{i_p' i_1'}) \\
    &=\E(X_{i_1 i_2} X_{i_1' i_2'}) \dots \E(X_{i_p i_1} X_{i_p' i_1'} X_{i_{k} i_{k+1}} X_{i_{k}' i_{k+1}'}) \\
    & \geq \gamma^{p-2}(3\gamma^2) \\
    &= 3 \gamma^{p}.
\end{align*}

The second equality comes from independence. If $p \geq 3$, then $X_{i_{1} i_{2}}$ is independent of $X_{i_{2} i_{3}}$, and in the case when $p = 2$, $X_{i_{1} i_{2}}$ is independent of $X_{i_{2} i_{1}}$ due to the lack of symmetry of circulant matrices. The first inequality comes from applying Wick's theorem to the expectation of the four Gaussians in the line above.

Thus following the rest of the proof of the generalized Toeplitz case completes the proof of the Lemma.

\done

\blem
\label{NormBoundC}
If $X_n$ is an $n \times n$ circulant matrix with correlated entries and $\lambda_n$ is the spectral norm of $X_n$, then $\E(\lambda_n) \leq C\sqrt{\log n}$, where $C$ is a constant independent of $n$. 
\elem

\prove
This is proved again via matrix concentration inequalities and a dilation argument. The corresponding $B_{n, k, l}$ matrices are defined by
\begin{equation*}\label{CBMatrix}
    (B_{n, k, l})_{ij} = \begin{cases}
        0 & j-i+1 \not \equiv k \mod n \\
        A_{i l} & j-i+1 \equiv k \mod n \\
    \end{cases}.
\end{equation*}

A similar computation to the reverse circulant case shows that 
\begin{equation*}
    (B_{n, k, l} B_{n, k, l}^T)_{ii} = A_{il}^2.
\end{equation*}
Since each of the rows of $A$ have length 1, $\sum_{l=1}^n\widetilde{B_{n,k,l}}^2 = I_{2n \times 2n}$ where $\widetilde{B_{n,k,l}}$ is the dilation
\begin{equation*}
   \widetilde{B_{n,k,l}} = \begin{bmatrix}
        0 & B_{n,k,l}^T \\
        B_{n,k,l} & 0 \
    \end{bmatrix}.
\end{equation*}
It then follows that $\left \Vert \sum_{k=1}^n \sum_{l=1}^n\widetilde{B_{n,k,l}}^2 \right \Vert^{1/2} \leq n^{1/2}$ and thus the lemma is proved.
\done

\subsection{Symmetric Circulant Matrices with Correlated Entries}

The following two Lemmas complete the proof. 

\blem
With $W_n$ as in Theorem \ref{MainTheorem} and $X_n$ a symmetric circulant matrix with correlated entries, $\Var(W_n) \geq Kn$ for some constant $K$ that only depends on $p$ and $\gamma$.
\elem

\prove
Note that when $|i - j| \leq \frac{n}{2}$, $X_{ij} = [\mathbf{X}_{|i-j|}]_{\min(i,j)}$ just as in the Toeplitz case. The argument used in Lemma \ref{VarBound} in the Toeplitz case only considers partitions where $|i_k - i_{k-1}| \leq \lceil \frac{n}{3p} \rceil$, so thus the argument applies in the case for symmetric circulant matrices with correlated entries, and we get $\Var(W_n) \geq Kn$ for the same constant $K$ as the Toeplitz case. 
\done

\blem
\label{NormBoundSC}
If $X_n$ is an $n \times n$ symmetric circulant matrix with correlated entries and $\lambda_n$ is the spectral norm of $X_n$, then $\E(\lambda_n) \leq C\sqrt{\log n}$, where $C$ is a constant independent of $n$. 
\elem

\prove
Note that for any symmetric circulant matrix $X_n$, we can write $X_n = Y_n + Z_n$ where $Y_n$ and $Z_n$ are Toeplitz matrices with correlated entries with diagonals 0 to $\lceil \frac{n}{2} \rceil$ replaced by 0's in $Z_n$ and diagonals $\lceil \frac{n}{2} \rceil + 1$ to $n-1$ replaced by 0's in $Y_n$. Then applying the same argument via matrix concentration inequalities as in the Toeplitz case to $Y_n$ and $Z_n$ yields
\begin{equation*}
    ||Y_n|| \leq C \sqrt{\log n} \qquad \text{and} \qquad ||Z_n|| \leq C \sqrt{\log n},
\end{equation*}
where $C$ is a constant independent of $n$. It then follows that $||X_n|| \leq 2C \sqrt{\log n}$ as desired.
\done

\subsection{Hankel Matrices with Correlated Entries}

We then have the following two Lemmas needed to complete the proof.

\blem
With $W_n$ as in Theorem \ref{MainTheorem} and $X_n$ a Hankel matrix with correlated entries, $\Var(W_n) \geq Kn$ for some constant $K$ that only depends on $p$ and $\gamma$.
\elem

\prove
 We can apply the same combinatorial argument as in the reverse circulant case here.
\done

\blem
\label{NormBoundH}
If $X_n$ is an $n \times n$ Hankel matrix with correlated entries and $\lambda_n$ is the spectral norm of $X_n$, then $\E(\lambda_n) \leq C\sqrt{\log n}$, where $C$ is a constant independent of $n$. 
\elem

\prove
Any Hankel matrix can be viewed as the the first $n \times n$ block of a $2n \times 2n$ reverse circulant matrix with correlated entries. Thus let $H_n$ be any Hankel matrix with correlated entries and $RC_{2n}$ be a reverse circulant matrix with its first $n \times n$ block being $H_n$. Then
\begin{equation*}
    ||H_n|| \leq ||RC_{2n}||.
\end{equation*}
Then by Lemma \ref{NormBoundRC},
\begin{equation*}
    \E||H_n|| \leq \E||RC_{2n}|| \leq C\sqrt{\log(2n)}  \leq \tilde C \sqrt{\log n}.
\end{equation*}
\done

\section{Fluctuations of Eigenvalues with Correlation Decay}
\label{DecaySection}

In this section we prove Theorem \ref{DecayTheorem}. The proof follows the same general structure as that of Theorem \ref{MainTheorem} by using Proposition \ref{ChatterjeeGaussian}. In order to prove the theorem as stated, we need tighter bounds on the operator norm. These bounds result from the sharper matrix concentration inequalities in \cite{BandeiraBoedihardjovanHandel2023}.

The set up is as follows. As in Theorem \ref{MatConcIneq}, let 
\begin{equation*}
    X = \sum_{i = 1}^{k} \gamma_{i} H_{i}
\end{equation*}
where $\gamma_{1}, \dots, \gamma_{k}$ are independent standard normal random variables, and $H_{1}, \dots, H_{k}$ are nonrandom Hermitian matrices of common dimension $n$, and again let 
\begin{equation*}
    \sigma(X) = ||\E(X^2)|| = \left \Vert \sum_{i = 1}^{k}H_{i}^{2} \right \Vert.
\end{equation*}
Further, define 
\begin{equation*}
    \nu(X) = ||\Cov(X)||,
\end{equation*}
where $\Cov(X)$ is viewed as a matrix in $M_{n^2 \times n^2}(\R)$ and $\Cov(X)_{ij, kl} = \Cov(X_{ij}, X_{kl})$.
We then have the following result, which is equation (1.11) in \cite{BandeiraBoedihardjovanHandel2023} (a corollary of their Theorem 1.2).
 
\bth 
\label{SharperMatConc}
With $X$, $\sigma(X)$, and $\nu(X)$ defined as above, $\E(||X||) \lesssim \sigma(X) + \nu(X)(\log n)^{3/2}$.
\eth

\begin{remark}
    Here $\nu(X)$  involves the correlation matrix for the random matrix $X$, i.e. it includes the $n^{-1/2}$ normalization for each of the entries. This is in contrast to Chatterjee's statement in Proposition \ref{ChatterjeeGaussian} in which we compute the operator norm of the covariance matrix of the Gaussian input sequence before we rescale by $n^{-1/2}$.
\end{remark}

Using the above theorem, we have the following Lemma which states that if we have any polynomial decay in the correlations, then the operator norm of any of the generalized models stays bounded in expectation as $n$ tends to infinity. 

\blem
\label{DecayNormBound}
Let $X_n$ be a random Toeplitz, circulant, reverse circulant, symmetric circulant, or Hankel matrix with correlated entries with standard Gaussian entries with covariances along each diagonal labeled by $c_k(i, j)$. Further, assume $c_k(i, j) = n^{-\alpha}$ for some $\alpha > 0$ for all $i, j, k$ with $i \neq j$. Then $\E(||X_n||)$ is bounded as $n \rightarrow \infty$.
\elem

\prove

In Theorem \ref{NormBoundIDID} and Lemmas \ref{NormBoundRC}, \ref{NormBoundC}, \ref{NormBoundSC}, and \ref{NormBoundH}, we showed that $\sigma(X)$ is bounded by a constant. In all of the models, each entry of $X_n$ is correlated with at most $2n-1$ other entries. Hence by the Gershgorin circle theorem, 
\begin{equation*}
    \nu(X) \leq \frac{1}{n} \left(1 + (2n-1)\sup_{k, i \neq j}c_k(i, j) \right).
\end{equation*}
The $n^{-1}$ factor comes from including the normalization factor $n^{-1/2}$ in $\Cov(X)$. Since $c_k(i, j) = o(n^{-\alpha})$, $\nu(x) = o(n^{-\alpha})$, and thus from Theorem \ref{SharperMatConc}, $\E(||X_n||) \leq C + o(n^{-\alpha})$.
\done

\begin{remark}
    Combining the above Lemma with the concentration arguments in Lemma \ref{OperatorNormMomentBound} shows that the operator norm of these generalized models has a sub-Gaussian distribution when there is any polynomial decay in the correlations among the entries.
\end{remark}

\prove (Proof of Theorem \ref{DecayTheorem})

The following proof works for any of the patterned matrix models with correlations considered in this paper. Let $X_n$ be a generalized random Toeplitz, circulant, reverse circulant, symmetric circulant, or Hankel matrix with correlated entries with standard Gaussian entries with covariances along each diagonal labeled by $c_k(i, j)$. Let $p$ be a positive integer. The structure of the proof is similar to that of the case where the correlations are bounded away from zero in that we show the fraction on the right hand side of equation (\ref{ChatEqn}) converges to zero. 

From the computations in Sections \ref{IDIDBound} and \ref{OtherModelSection}, we know that the expectation of the operator norm of $X_n$ is at most $\sqrt{\log n}$. Thus from Lemma \ref{OperatorNormMomentBound} and Lemma \ref{DecayNormBound}, the term $ab$ in equation (\ref{ChatEqn}) is bounded by $p^{p+3}C^p$ where $C$ is independent of $n$. 

Next we bound $\sigma^2$ from below by a constant. Recall from the proof of Lemma \ref{VarBound} that the covariance of arbitrary products of the entries of $X_n$ is nonnegative. Thus
\begin{align*}
    \Var(W_n) &= n^{-p}\Var \left(\sum_{1 \leq i_1, \dots, i_p \leq n} X_{i_1i_2}X_{i_2i_3}\dots X_{i_pi_1}\right) \\
    &\geq n^{-p}\sum_{1 \leq i_1, \dots, i_p \leq n} \Var(X_{i_1i_2}X_{i_2i_3}\dots X_{i_pi_1}).
\end{align*}
In order to apply a partition argument via Wick's Theorem to the variances in the sum, it will be again easier to relabel the $2p$ random variables $X_{i_1i_2},X_{i_2i_3},\dots X_{i_pi_1}, X_{i_1i_2},X_{i_2i_3},\dots, X_{i_pi_1}$ by enumerating them from $1$ to $2p$. Further, define the partition $\tau \coloneqq \{\{1,p+1\}, \{2, p+2\}, \dots, \{p, 2p\}\}$. Then we have
\begin{align*}
    \Var(X_{i_1i_2}X_{i_2i_3}\dots X_{i_pi_1}) &= \sum_{\pi \in P_2(2p)} \prod_{\{i, j\} \in \pi} \E(X_i X_j) \\
    & - \left(\sum_{\pi \in P_2(p)} \prod_{\{i, j\} \in \pi} \E(X_i X_j) \right) \left(\sum_{\pi \in P_2(p)} \prod_{\{i, j\} \in \pi} \E(X_{i+p} X_{j+p}) \right). \\
\end{align*}
Now note that every term in the double sum on the right appears in the double sum on the left (when $p$ is odd the terms on the right are 0). However, $\tau$ does not appear in the double sum due its blocks crossing between the sets $\{1, \dots, p\}$ and $\{p+1, \dots, 2p\}$. Hence
\begin{align*}
     \Var(X_{i_1i_2}X_{i_2i_3}\dots X_{i_pi_1}) \geq \prod_{\{i, j\} \in \tau} \E(X_i X_j) = 1,
\end{align*}
and it follows that 
\begin{equation*}
    \sigma^2 = \Var(W_n) \geq 1.
\end{equation*}

The last part of the proof is to bound the operator norm of the covariance matrix $\Sigma$. Following the proof of Lemma \ref{CovarianceMatrixBoundOtherModels},
\begin{equation*}
    ||\Sigma|| \leq 1 + (2n-1) \sup_{k, i \neq j} c_k(i, j).
\end{equation*}

Under the assumption that $c_k(i, j) = o(n^{-1/3})$, $||\Sigma|| = o(n^{2/3})$ and $||\Sigma||^{3/2} = o(n)$. Then combining this with the lower bound for the variance and upper bound on $ab$ via Lemmas \ref{DecayNormBound} and \ref{OperatorNormMomentBound}, and plugging into equation (\ref{ChatEqn}) we get
\begin{equation*}
    d_{TV}(W_n, Z_n) = o(1).
\end{equation*}

\done
\section{Towards Universality}
\label{UniversalitySection}

\subsection{Removing the Conditions on the Covariances}

One obvious extension of Theorem \ref{MainTheorem} and Theorem \ref{DecayTheorem} is to allow for completely arbitrary correlation structures. Hence we have the following conjecture.

\begin{conj}
\label{ConbjectureArbitraryCorrelation}
Let $X_n$ be an $n \times n$ generalized random Toeplitz, circulant or symmetric circulant matrix with standard Gaussian entries. Let $p \geq 2$ be a positive integer and let $W_n = \mathrm{Tr}(X_n^p)$. Then as $ n \rightarrow \infty$,
\begin{equation*}
\frac{W_n - \E (W_n)}{\sqrt{\Var(W_n)}} \text{ converges in total variation to } N(0,1).
\end{equation*}

In the case when $X_n$ is a reverse circulant or Hankel matrix with correlated entries, the theorem still holds under the further assumption that $p$ is restricted to be an even positive integer.
\end{conj}

In the proofs of Theorems \ref{MainTheorem} and \ref{DecayTheorem} the fraction in equation (\ref{ChatEqn}) vanishes as $n \rightarrow \infty$ for different reasons. Hence in order to prove Conjecture \ref{ConbjectureArbitraryCorrelation}, one would need to adopt new methods to more precisely simultaneously bound all three terms that we bounded in our proofs, or an entirely new approach is needed. Furthermore, there is the fact that Theorem \ref{DecayTheorem} holds for Hankel matrices with correlated entries when $p$ is odd, but the linear eigenvalue statistics of a Hankel matrix for odd $p$ converge to a non-Gaussian limit (see Theorem 2 of \cite{KumarMauryaSaha2022}). Hence there is most likely a phase transition for the limiting statistics of odd monomial test functions for linear eigenvalue statistics of Hankel matrices, and a similar phenomenon may occur in the case of reverse circulant matrices. We also did not consider generalized models with negative correlations.

\subsection{Moving Beyond Gaussian Entries}

In \cite{AdhikariSaha2018}, Adhikari and Saha were able to prove a Gaussian central limit theorem for circulant matrices with sub-Gaussian entries. They further assumed that the entries belonged to a certain class of random variables whose laws can be written as $C^2$ functions with bounded derivatives of a standard Gaussian, denoted $\mathcal{L}(c_1, c_2)$ where $c_1$ and $c_2$ are the respective bounds for the first and second derivatives (see  Definition 2.1 of \cite{Chatterjee2009}). There are multiple equivalent definitions of sub-Gaussian random variables (see Proposition 2.5.2 of \cite{Vershynin2018}), and one such definition is the following.

\bde
A random variable $X$ is said to be sub-Gaussian if there exists a constant $K$ such that the tails of x satisfy
\begin{equation*}
    \P(|X| \geq t) \leq 2e^{-t^2/K^2} \quad \text{for all $t \geq 0$}.
\end{equation*}
\ede

We conjecture that similar extensions to generalized patterned random models studied with sub-Gaussian entries hold due to the sub-Gaussian tails inducing similar $\sqrt{\log n}$ bounds on the operator norm. Meckes proved in \cite{Meckes2007} that the operator norm of a random Toeplitz matrix (and other models considered in this paper) with sub-Gaussian entries is asymptotically bounded above by $\sqrt{\log n}$. Due to concentration of measure, the original models (with ungeneralized correlation structure) have the most dependence among their entries and thus are the ``worst-case" scenario for bounding the operator norm. Thus it is natural to believe that for any of the generalized models in this paper with sub-Gaussian entries, the operator norm is asymptotically bounded above by $\sqrt{\log n}$. We then have the following conjecture.
\begin{conj}
Let $X_n$ be an $n \times n$ generalized random Toeplitz, circulant or symmetric circulant matrix with symmetric standard (i.e., mean zero and variance one) sub-Gaussian entries. Let $p \geq 2$ be a positive integer and let $W_n = \mathrm{Tr}(X_n^p)$. Then as $ n \rightarrow \infty$,
\begin{equation*}
\frac{W_n - \E (W_n)}{\sqrt{\Var(W_n)}} \text{ converges in total variation to } N(0,1).
\end{equation*}

In the case when $X_n$ is a generalized random reverse circulant or Hankel matrix, the theorem still holds under the further assumption that $p$ is restricted to be an even positive integer.
\end{conj}

In order to prove this, new methods will most likely need to be used. Proposition \ref{ChatterjeeGaussian} from \cite{Chatterjee2009} only holds for Gaussian entries, though Theorems 2.2 and 3.1 of the same paper can be applied to more general distributions of the entries (and this was how Adhikari and Saha proved their sub-Gaussian universality result in \cite{AdhikariSaha2018}). However, these results only hold under the assumption that each entry of the matrix can be written as a function of independent random variables. It is not straightforward to construct a random matrix with a general correlation structure and sub-Gaussian entries in this manner. 

\printbibliography

\end{document}